\documentclass{article}

\usepackage[top=1.5in,left=1.25in,right=1.25in,bottom=1.5in]{geometry}

\usepackage{latexsym}
\usepackage[pdftex]{graphicx}
\usepackage{epsfig, ecltree, epic, eepic}
\usepackage{enumerate,amsmath,amssymb,dsfont,pifont,amsthm}
\usepackage{xcolor}
\usepackage{ifthen}
\usepackage{graphicx}
\usepackage[arrow, matrix, curve]{xy}
\usepackage{amsthm,amsmath}
\usepackage{stmaryrd}
\usepackage{comment}
\usepackage{enumerate,amssymb,dsfont,pifont}
\usepackage{tikz}
\usepackage{hyperref}
\hypersetup{
    colorlinks,
    citecolor=black,
    filecolor=black,
    linkcolor=black,
    urlcolor=black
}

\theoremstyle{plain}
\newtheorem{theorem}{Theorem}[section]

\newtheorem{lemma}[theorem]{Lemma}
\newtheorem{proposition}[theorem]{Proposition}

\theoremstyle{definition}
\newtheorem{definition}[theorem]{Definition}

\theoremstyle{remark}

\newtheorem{example}[theorem]{Example}

\newcommand{\bbc}{\mathbb{C}}
\newcommand{\bbr}{\mathbb{R}}

\newcommand{\bbp}{\mathbb{P}}
\newcommand{\bbe}{\mathbb{E}}

\newcommand{\bbn}{\mathbb{N}}

\newcommand{\cm}{\mathcal{M}}

\newcommand{\abs}[1]{\left| #1 \right|}
\newcommand{\norm}[1]{\left\| #1 \right\|}

\newcommand{\roi}[2]{\left[ #1 , #2 \right[}

\makeatletter
\newcommand*\Si{\mathpalette\Si@{.5}}
\newcommand*\Si@[2]{\mathbin{\vcenter{\hbox{\scalebox{#2}{$\m@th#1\bullet$}}}}}
\makeatother

\newcommand{\ls}{\llbracket}
\newcommand{\rs}{\rrbracket}

\DeclareMathOperator{\sign}{sign}

\renewcommand{\leq}{\leqslant}
\renewcommand{\geq}{\geqslant}
\renewcommand{\Re}{\ensuremath{{\rm Re\,}}}
\renewcommand{\Im}{\ensuremath{{\rm Im\,}}}

\begin{document}

\newcounter{farbe}
    \setcounter{farbe}{0} 

\ifthenelse{\value{farbe}=1}{\newcommand{\cor}{\color{blue}}}{\newcommand{\cor}{\color{black}}}

\allowdisplaybreaks

\title{\bfseries On the r\^ole of singular functions in extending the probabilistic symbol to its most general class}

\author{%
    \textsc{Sebastian Rickelhoff}%
    \thanks{University of Siegen, Department of Mathematics, Walter-Flex-Straße 3,
              D-57072 Siegen, Germany,
              \texttt{Sebastian.Rickelhoff@uni-siegen.de}, 
              phone: +49 (271) 740-3575.} \ and
                  \textsc{Alexander Schnurr}%
    \thanks{University of Siegen, Department of Mathematics,Walter-Flex-Straße 3,
               D-57072 Siegen, Germany,
              \texttt{Schnurr@mathematik.uni-siegen.de}, 
              phone: +49 (271) 740-3806.}
    }

\date{\today}

\maketitle
\begin{abstract} The probabilistic symbol is the right-hand side derivative of the characteristic functions corresponding to the one-dimensional marginals of a stochastic process. This object, as long as the derivative exists, provides crucial information concerning the stochastic process. For a L\'evy process, one obtains the characteristic exponent while the symbol of a (rich) Feller process coincides with the classical symbol which is well known from the theory of pseudodifferential operators. \\
Leaving these classes behind, the most general class of processes for which the symbol still exists are Lévy-type processes. 
It has been an open question, whether further generalizations are possible within the framework of Markov processes.  We answer this question in the present article: within the class of Hunt semimartingales, L\'evy-type processes are exactly those for which the probabilistic symbol exists. Leaving {\cor Hunt} behind, one can construct processes admitting a symbol. However, we show, that the applicability of the symbol might be lost for these processes.  Surprisingly, in our proofs the upper and lower Dini derivatives corresponding to certain singular functions play an important r\^ole. 
\end{abstract}



\section{Introduction and History} \label{sec:intro}

Consider a stochastic process $X=(X_t)_{t\geq 0}$ starting in $x\in\bbr^d$. The \emph{probabilistic symbol} is the right-hand side derivative of the characteristic function of $(X_t-x)_{t\geq 0}$, that is, 
\begin{align} \label{probsymbol1} 
     q(x,\xi):=- \lim_{t\downarrow 0} \frac{\bbe^xe^{i(X_t-x)'\xi}-1}{t}.
\end{align}
The characteristic function carries all the information of the marginal distribution \emph{at time} $t$. Hence, it is not surprising, that the right-hand side derivative \eqref{probsymbol1} might be used to analyze the development of the process \emph{over time.}  From the point of view of Markov processes, Formula \eqref{probsymbol1} can be interpreted as follows: one plugs a (complex) exponential function into the generator of the process.  
In this spirit the formula was used for the first time in 1998 by Jacob \cite{nielsursprung} in order to
calculate the functional analytic symbol of a Feller process satisfying some additional properties (see below). 
Since then, this formula, related maximal inequalities and results on path properties have been generalized to further classes of processes. 
{\cor All of the classes mentioned here, will be properly defined further down.}
Within the Markovian framework the most general class for which this symbol is known to exist is the class of L\'evy-type processes \cite{levymatters3} which are sometimes called It\^o processes in the sense of Çinlar, Jacod, Protter and Sharpe \cite{vierleute}.  
A canonical class for a further generalization would be Hunt semimartingales. In \cite{withMartynas} we have shown that within this class, there exist processes which do not admit a symbol. 
It has been an open question for over 10 years, whether or not L\'evy-type processes are the most general class for which the symbol exists and can be used. This is also point (d) in the open problems section of the monograph \cite{levymatters3}.  We answer this question in the present paper: we prove the existence of a Markov semimartingale, which is not a L\'evy-type process but which admits a symbol. This process gives rise to a class of processes of this kind. 
Furthermore, we show that there is no Hunt process admitting a symbol which is not a L\'evy-type process (Theorem \ref{thm:main}). In this sense the above conjecture is true.  

Summing up: Under very weak conditions every L\'evy-type process admits a symbol. There are no other Hunt semimartingales admitting a symbol. But if we move on another step and consider Markov semimartingales we again find processes for which the symbol does exist. However, in this case the symbol does not carry the valuable information on the process anymore. 

Let us first fix some notation: {\cor If not mentioned otherwise every stochastic process $X:=(X_t)_{t \geq 0}$ on a stochastic basis $(\Omega, \mathcal{A}, (\mathcal{A}_t)_{t \geq 0},\bbp)$ takes values in $(\bbr^d, \mathcal{B}(\bbr^d))$ and is c\'{a}dl\'{a}g. Here, $\mathcal{B}(\bbr^d)$ is the $\sigma$-field of Lebesgue sets.  We denote by $\mathcal{F}^0_t:= \sigma( X_s: s \leq t)$ and $\mathcal{F}^0:= \sigma(X_s: s\geq 0)$ the \textit{natural} filtration of $X$.  Moreover, we call $\Delta X_t := X_t - \lim_{s \uparrow t} X_s$ the jump of the process at time $t\geq 0$, and for a stopping time $\tau$ we call $X^\tau:= X1_{\ls 0, \tau \rs}+ X_\tau 1_{\ls \tau, \infty \ls}$ the stopped process. Here,  the stochastic interval $\ls \tau , \sigma \ls$ for two stopping times $\tau, \sigma$ is defined by $\{ (\omega, t) \in \Omega\times \bbr_+: \tau(\omega) \leq t < \sigma(\omega)\}$. The stochastic intervals $\ls \tau , \sigma \rs$, $\rs \tau , \sigma \ls$, $\rs \tau , \sigma \rs$ are defined alike.} \\
A {\cor (strong)} \emph{Markov process} $(\Omega, \cm, (\cm_t)_{t\geq 0}, (X_t)_{t\geq 0}, (\theta_t)_{t \geq 0},\bbp^x)_{x\in\bbr^d}$ is defined in the sense of Blumen-thal-Getoor (cf. \cite{blumenthalget}).{\cor We assume that $\mathcal{M}= \overline{\mathcal{M}}$ and $(\mathcal{M}_t)_{t \geq 0}= (\overline{\mathcal{M}}_t)_{t \geq 0}$ where $\overline{\mathcal{M}}$ resp. $\overline{\mathcal{M}}_t$ is the completion of $\mathcal{M}$ resp.  $\mathcal{M}_t$ with respect to $\{ \bbp^x: x \in \bbr^d\}$ (cf. \cite{blumenthalget} Section I.5). Moreover, we assume all filtrations  encountered in the following to be right-continuous}. Being Markovian means in particular that the following formula holds
\begin{align} \label{universalmp}
  P_{s,t}^w(x,B) = P_{0,t-s}^{x}(x,B)=:P_{t-s}(x,B)
\end{align}
where $P_{s,t}^w(x,B)$ is the regular version of $\bbp^w(X_t\in B \, | \, X_s=x)$ with $w,x\in\bbr^d$, $s\leq t$ and $B$ is a Borel set in $\bbr^d$. Processes of this kind are sometimes called Markov families or {\cor (strong) Markov process with transition kernel} (cf. \cite{bauerwt}, \cite{niels3}) in the literature.  {\cor To avoid possible misunderstandings concerning this notion we only use the term (strong) Markov process subsequently}. We assume that every Markov process is normal, that is,  $\bbp^x(X_0=x)=1$.
We associate a semigroup $(T_t)_{t\geq 0}$ of operators on $B_b(\bbr^d)$ with every Markov process by setting
\[
    T_t u(x):= \bbe^x u(X_t), \quad t\geq 0,\; x\in \bbr^d.
\]
Let us denote by $C_\infty(\bbr^d)$ the space of all functions $u:\bbr^d\to\bbr$ which are continuous and vanish at
infinity. We call $(T_t)_{t\geq 0}$ a Feller
semigroup and $(X_t)_{t\geq 0}$ a \emph{Feller process}, if the semigroup is strongly continuous, that is, $\norm{T_tu-u}_\infty\to 0$ for $t\downarrow 0$, and if the following condition is satisfied:
\begin{align}
  T_t:C_\infty(\bbr^d) \to C_\infty(\bbr^d) \text{ for every }t\geq 0.
\end{align}

The \emph{generator} of the Feller process $(A,D(A))$ is the closed operator given by
\begin{gather}\label{generator}
    Au:=\lim_{t \downarrow 0} \frac{T_t u -u}{t} \qquad\text{for\ \ } u\in D(A)
\end{gather}
where the \emph{domain} $D(A)$ consists of all $u\in C_\infty(\bbr^d)$ for which the limit \eqref{generator} exists uniformly. A Feller process is called \emph{rich} if $C_c^\infty(\bbr^d)\subseteq D(A)$, that is, the test functions are contained in the domain of the generator. 
{\cor The generator of those processes fulfills the so called positive maximum principle N. Jabob came up with the idea to use this fact combined with a classical result due to P.\ Courr\`ege \cite{courrege} (cf. in this context also \cite{vw1} and \cite{vw2}) to derive the following}: the generator of a rich Feller process is (restricted to the test functions $C_c^\infty(\bbr^d)$) a \emph{pseudo differential operator}, i.e.,\ $A$ can be written as
\begin{gather} \label{pseudo}
    Au(x)= - \int_{\bbr^d} e^{ix'\xi} q(x,\xi) \widehat{u}(\xi) \,d\xi, \qquad u\in C_c^\infty(\bbr^d)
\end{gather}
where $\widehat{u}(\xi)=(2\pi)^{-d}\int e^{-iy'\xi}u(y) dy$ denotes the Fourier transform. The function $q:\bbr^d \times \bbr^d \to \bbc$ is locally bounded and, for fixed $x$, a continuous negative definite function in the sense of Schoenberg in the co-variable $\xi$ (cf. \cite{bergforst} Chapter II {\cor and \cite{niels1} Section 3.7}). This is equivalent to the fact that it admits a L\'evy-Khintchine representation
\begin{align} \label{lkfx}
  q(x,\xi)=
  -i \ell(x)'  \xi + \frac{1}{2} \xi'Q(x) \xi -\int_{w\neq 0} \left( e^{i \xi' w} 
  -1 - i \xi' w \cdot \chi(w)\right)N(x,dw)
\end{align}
where $\ell(x)=(\ell^{(j)}(x))_{1\leq j \leq d} \in \bbr^d$, $Q(x)=(Q^{jk}(x))_{1\leq j,k \leq d}$ is a symmetric positive semidefinite matrix and $N(x,dw)$ is a measure on $\bbr^d\setminus\{0\}$ such that $\int_{w\neq 0} (1 \wedge \norm{w}^2) \,N(x,dw) < \infty$ and $\chi:\bbr^d\to\bbr$ is a cut-off function, {\cor i.e.  $\chi$ is measurable with compact support and equal to 1 in a neighborhood of zero}. The function $q:\bbr^d\times \bbr^d\to \bbc$ which is often written as $q(x,\xi)$ is called the \emph{symbol} of the operator. For details on the rich theory of the interplay between processes and their symbols we refer the reader to \cite{levymatters3}  and \cite{niels3}.

In \cite{nielsursprung} and for a more general class of Feller processes in \cite{schilling98pos} it was shown that the formula \eqref{probsymbol1} can be used to calculate the symbol $q(x,\xi)$ of the process.  For L\'evy processes this is the characteristic exponent. 

{\cor We call a strong Markov process \textit{Hunt process} if it is quasi left-continuous (cf. \cite{jacodshir} Definition I.2.25). }
A Markov process $X$ is called \emph{Markov semimartingale}, if $X$ is for every $\bbp^x$ a semimartingale. If such a process is in addition a Hunt process, it is called \emph{Hunt semimartingale}. In \cite{mydiss} we have shown that every rich Feller process is a Hunt semimartingale and even a \emph{L\'evy-type process}, that is, a Markov semimartingale with characteristics of the following form
\begin{align} \begin{split} \label{chars}
  B^{(j)}_t(\omega)&=\int_0^t \ell^{(j)}(X_s(\omega)) \ ds \\
  C^{jk}_t(\omega)&=\int_0^t Q^{jk} (X_s(\omega)) \ ds \\
  \nu(\omega;ds,dy)&=N(X_s(\omega),dy)\ ds
\end{split}\end{align}
for every starting point $x\in\bbr^d$ with respect to a fixed cut-off function $\chi$. Here $\ell$, $Q$ and $N$ are as above. The triplet $(\ell,Q,n:=\int_{y\neq 0} (1\wedge \norm{y}^2) \ N(\cdot,dy)$) is called the \emph{differential characteristics} of the process. The class of semimartingales which are not necessarily Markovian, but having characteristics as described in \eqref{chars} are called \emph{homogeneous diffusions with jumps}  (cf. \cite{jacodshir} Definition III.2.18). 

Since the symbol and some related {\cor indices} have proven to be useful in analyzing global and path properties, we have generalized the notion of the symbol using a generalization of formula \eqref{probsymbol1}, which can be used even for processes which do not admit a generator, namely homogeneous diffusions with jumps.

\begin{definition} \label{def:symbol}
Let $X$ be a stochastic process. Fix a starting point $x$ and define $\sigma=\sigma^x_k$ to be the first exit time from a compact ball with radius $k$ and $x$ as center
  \[
    \sigma:= \inf \{ t \geq 0 : \| X_t-x \| > k \}.
  \]
For $\xi \in \bbr^d$ we call $p:\bbr^d\times \bbr^d\to \bbc$ given by
\begin{align} \label{probsymbol2} 
     p(x,\xi):=- \lim_{t\downarrow 0}\bbe^x \frac{e^{i(X^\sigma_t-x)'\xi}-1}{t}
\end{align}
the \emph{(probabilistic) symbol} of the process, if the limit exists and coincides for every choice of $k>0$.
\end{definition}

In \cite{generalizedindizes} we have shown that the symbol exists for homogeneous diffusions with jumps having differential characteristics which are locally bounded and finely continuous (cf. \cite{blumenthalget},  Theorem II.4.8). Like the functional analytic symbol in the Feller case it is a state-space dependent continuous negative definite function in the sense of Schoenberg. Instead of the compact balls, one can use general compact neighborhoods. In the definition, it is demanded that the limit exists for all $k>0$. If this is the case, it is easily seen that the limit is indeed the same for every $k>0$ and even for every compact neighborhood of $x$. The stopping time is needed in order to have the process bounded at least on $\ls 0, \sigma \ls $.
In order to calculate the symbol, \emph{fine continuity} (and local boundedness) are the most general requirements. However, since classical continuity is much more natural, we demand this property from now on. Continuity implies both: fine continuity and local boundedness.  Hence, Theorem 3.6 of \cite{generalizedindizes} yields that the symbol exists for every L\'evy-type process with continuous differential characteristics. 


For the readers convenience we recall some concepts and results on advanced calculus on the real line (cf.  \cite{kankrue} and \cite{natanson}).

{\cor
\begin{definition}\label{def:singular}
A function $f:\bbr \to \bbr$ is called singular if it is non-constant and whose derivative exists and is zero almost everywhere.
\end{definition} 
Almost everywhere is meant w.r.t. the Lebesgue sets, that is, we are working on the completion of the Borel sets w.r.t. Lebesgue measure.
}

\begin{definition}\label{def:dinider}
Let $f:[a,b] \to \bbr$ be a function, and let $x_0 \in [a,b]$ for $a,b \in \bbr$ and $a<b$. Then 
\begin{itemize}
\item[(a.)] the \emph{upper {\cor resp. lower} right Dini derivative} of $f$ at $x_0$ is defined by 
$$
D^+f(t):= \limsup_{x \downarrow x_0} \frac{f(x)-f(x_0)}{x-x_0} \text{ {\cor resp. }} D_+f(t):= \liminf_{x \downarrow x_0} \frac{f(x)-f(x_0)}{x-x_0},
$$
\item[(b.)] the \emph{upper {\cor resp. lower} left Dini derivative} of $f$ at $x_0$ is defined by 
$$
D^-f(t):= \limsup_{x \uparrow x_0} \frac{f(x)-f(x_0)}{x-x_0} \text{{ \cor resp. }} D_-f(t):= \liminf_{x \uparrow x_0} \frac{f(x)-f(x_0)}{x-x_0}.
$$
\end{itemize}
\end{definition}

The concept of singular functions (Definition 7.1.44 \cite{kankrue}) is fundamental for the reasoning in the subsequent sections. More specifically, in order to prove the main theorem of this paper,  we use the fact that a singular function possesses at least one point where a Dini derivative is infinite.  Although, this seems to be clear since the most prominent examples like the Cantor function or Minkowski's question mark function on $[0,1]$ possess infinite derivatives in \textit{all} points where the derivative is not zero, it was shown in \cite{spanier} that there exists a singular function with a derivative that takes non-zero finite values on an uncountable dense set whose intersection with any interval $(a,b)$ possesses Hausdorff dimension one.

\begin{proposition} \label{prop:derived}
Let $g:\bbr_+ \to \bbr$ be an increasing, singular function which is not constantly zero. Then there exist at least one point where $g$ possesses an infinite upper right (resp. upper left, lower right, lower left) Dini derivative.
\end{proposition}

{\cor In order to prove this proposition one can combine two results which can be found in \cite{saks}, namely Theorem (4.6) in Chapter IX and Theorem (6.7) in Chapter VII. We include a different proof in the appendix, to make the present article more self contained. }

\section{Main Result} \label{sec:general}

A Hunt semimartingale is a strong Markov process $(\Omega, \mathcal{M}, (\mathcal{M}_t)_{t \geq 0}, (X_t)_{t\geq 0}, (\theta_t)_{t \geq 0}, \bbp^x)_{x \in \bbr^d}$ which is also a $\bbp^x$-semimartingale and quasi-left-continuous.  
Within this framework, the following definition is useful.

\begin{definition}\label{def:addi}
Let $(\Omega, \mathcal{M}, (\mathcal{M}_t)_{t \geq 0}, (X_t)_{t\geq 0}, (\theta_t)_{t \geq 0}, \bbp^x)_{x \in \bbr^d}$ be a Markov process. 
\begin{itemize} 
\item[(a.)] We call a process $Y=(Y_t)_{t\geq 0}$ {\cor \emph{additive functional (AF)}} if 
\begin{itemize}
\item[(i.)] $Y_0=0$  {\cor $\bbp^x$-a.s.  for all $x \in \bbr^d$.}
\item[(ii.)] for every $s,t \geq 0$ we have 
$$
Y_{s+t} = Y_s+Y_t \circ \theta_s \quad {\cor \bbp^x\text{-a.s. for all }x \in \bbr^d.}
$$
\end{itemize}
\item[(b.)] We call $Y$ \emph{strongly {\cor additive functional}} if instead of $(ii.)$ 
\begin{itemize}
\item[(ii.')] for every $t \geq 0$ and all $(\mathcal{M}_t)$-stopping times $\rho$
$$
Y_{\rho+t}=Y_{\rho}+ Y_t \circ \theta_\rho 
$$
holds true {\cor $\bbp^x$-a.s.  for all $x \in \bbr^d$.}
\end{itemize}
{\cor
\item[(c.)] We call an $AF$ $Y$ \textit{perfect} if 
$$
\bigcup_{s,t \geq 0} \left\{Y_{s+t} \neq Y_t + Y_s \circ \theta_t \right\} 
$$
is null for all $\bbp^x, x \in \bbr^d$. }
\end{itemize}
\end{definition}
\noindent
We denote by $\mathcal{V}$ (resp. $\mathcal{V}^+$) the set of all real-valued processes that are càdlàg, adapted and possess paths of finite variation (resp. increasing paths), and we denote by $\mathcal{V}_{ad}$ (resp. $\mathcal{V}^+_{ad}$) the set of all additive functionals belonging to $\mathcal{V}$ (resp.  $\mathcal{V}^+$).

By Theorem 6.27 of \cite{vierleute} there exists for every Hunt semimartingale a continuous, strictly increasing, strongly {\cor additive functional} $F \in \mathcal{V}$ (Although,  in 6.27 $F$ is stated to be additive only, by the remark below 6.14 , one can assume $F$ to be strongly additive), measurable functions $\ell: \bbr^d \to \bbr^d$, $Q:\bbr^d \to \bbr^{d \times d}$ and a positive kernel $N$ from $(\bbr^d, \mathcal{B}(\bbr^d))$ into $(\bbr^d, \mathcal{B}(\bbr^d))$ such that the characteristics of $X$ are of the form 
\begin{align}
B^j_t&= \int_0^t \ell^{(j)}(X_s)  \ dF_s \label{b} \\
C^{jk}_t&= \int_0^t Q^{jk}(X_s)  \ dF_s  \label{c} \\
\nu(\omega; dt,dx) &= dF_t(\omega) N(X_t(\omega),dx).\label{K}
\end{align}

The proof of the following lemma is very similar to the one of Theorem 4.4 of \cite{mydiss}.  Of course it is possible to generalize the statement to the multidimensional case but we omit that since it is notationally more involved.
\begin{lemma}\label{Lemma 1}
Let $(\Omega, \mathcal{M}, (\mathcal{M}_t)_{t \geq 0}, (X_t)_{t\geq 0}, (\theta_t)_{t \geq 0}, \bbp^x)_{x \in \bbr}$ be a Hunt semimartingale with characteristics stated in (\ref{b})-(\ref{K}), and let 
$$
\sigma:= \inf \{ t \geq 0 : \| X_t-x \| > k \}
$$
for some starting point $x$ and $k \geq 0$.
Then the following equation holds true for $\xi \in \bbr$:
\begin{align*}
 \bbe^x \left( e^{i(X^\sigma_t-x)\xi } -1 \right)
 &=\bbe^x \left[ \int_{0}^t 1_{\ls 0,\sigma \ls}\left(i \xi e^{i(X^\sigma_{s-}-x)\xi}\ell(X_s)-  \frac{1}{2}\xi^2 e^{i(X^\sigma_{s-}-x)\xi} Q(X_s) \right.\right.\\
&+ \left. \left.\int_{ \{y\neq 0\}}  \Big(e^{i (X_{s-}-x)\xi}(e^{i\xi'y}-1 -i \xi'y\chi(y))\Big) \ N(X_s(\omega),dy) \right)\ dF_s(\omega)\right].
\end{align*}
\end{lemma}
\begin{proof}
The left-continuous process $X^\sigma_{t-}$ is bounded, the stopped jumps $(\Delta X)^\sigma$ are the jumps of the stopped process $(\Delta X^\sigma)$ and $X^\sigma$ admits the stopped characteristics:
\begin{align*}
B^\sigma_t(\omega)& = \int_0^t \ell (X_s(\omega)) 1_{\ls 0, \sigma \rs}(\omega, s) \ dF_s(\omega) \\
C_t^\sigma(\omega)&=\int_0^t Q(X_s(\omega)) 1_{\ls 0, \sigma \rs}(\omega, s) \ dF_s(\omega) \\
\nu^\sigma(\omega;dF_s,dy)&:=1_{\ls 0, \sigma \rs}(\omega, s) \ N(X_s(\omega),dy) \ dF_s(\omega).
\end{align*}
One can now set the integrand at the right endpoint of the stochastic support to zero, as we are integrating with respect to a continuous measure:
\begin{align*}
{\cor B^\sigma_t}(\omega)&=\int_0^t \ell (X_s(\omega)) 1_{\ls 0, \sigma \ls}(\omega, s) \ dF_s(\omega) \\
C_t^\sigma(\omega)&=\int_0^t Q(X_s(\omega)) 1_{\ls 0, \sigma \ls}(\omega, s) \ dF_s(\omega) \\
\nu^\sigma(\omega;dF_s,dy)&=1_{\ls 0, \sigma \ls}(\omega, s) \ N(X_s(\omega),dy) \ dF_s(\omega).
\end{align*}
Using It\^o's formula we obtain
\begin{align*}
 \bbe^x \left( e^{i(X^\sigma_t-x)'\xi } -1 \right)
&= \bbe^x \left(\int_{0+}^t i \xi e^{i(X^\sigma_{s-}-x)'\xi} \ dX^\sigma_s \right)\\
&+  \bbe^x\left(\frac{1}{2} \int_{0+}^t-\xi^2 e^{i(X^\sigma_{s-}-x)'\xi} \ d[X^\sigma,X^\sigma]_s^c \right)\\
&+  \bbe^x\left(e^{-ix'\xi} \sum_{0<s\leq t} \Big(e^{i \xi'X_s^\sigma}-e^{i\xi'X^\sigma_{s-}}-i\xi e^{i\xi' X^\sigma_{s-}} \Delta X^\sigma_s \Big)\right)\\
&=  \bbe^x \left(\int_{0+}^t i \xi e^{i(X^\sigma_{s-}-x)'\xi} \ dX^\sigma_s \right)\\
&+  \bbe^x\left(\frac{1}{2} \int_{0+}^t-\xi^2 e^{i(X^\sigma_{s-}-x)'\xi} \ dC^\sigma_s \right)\\
&+ \bbe^x\left( \int_{]0,t]\times \{y\neq 0\}}  \Big(e^{i (X_{s-}-x)'\xi}(e^{i\xi'y}-1 -i \xi'y\chi(y))\Big)\ \mu^{X^\sigma}(\cdot;ds,dy)\right) \\
&+ \bbe^x\left( \int_{]0,t]\times \{y\neq 0\}} \Big(  e^{i (X_{s-}-x)'\xi}( -i\xi'y \cdot(1-\chi(y)))\Big) \ \mu^{X^\sigma}(\cdot;ds,dy) \right).
\end{align*}
Now, we want to consider the first summand from above. To this end,  we use the canonical decomposition of the semimartingale (see \cite{jacodshir} Theorem II.2.34):
$$
X_t=X_0 + \left(X^\sigma_t\right)^c + \int_0^t \chi(y)y \ \left(\mu^{X^\sigma}(\cdot \ ;ds,dy)-\nu^\sigma(\cdot \ ;ds,dy)\right) +\int_0^t y(1-\chi(y)) \ \mu^{X^\sigma}(\cdot \ ;ds,dy) +B^\sigma_t.
$$
Thus, the linearity of the integral provides
\begin{align*}
\bbe^x \left(\int_{0+}^t i \xi e^{i(X^\sigma_{s-}-x)'\xi} \ dX^\sigma_s \right)&= \bbe^x \left(\int_{0+}^t i \xi e^{i(X^\sigma_{s-}-x)'\xi} \ dX^{\sigma,c}_s \right)\\
&+ \bbe^x\left( \int_0^t\int_{y\neq 0} i\xi \Big(e^{i(X_{s-}-x)'\xi}  \chi(y)y \Big) \mu^{X^\sigma}(\cdot;ds,dy)-\nu^\sigma(\cdot;ds,dy))\right) \\
&+ \bbe^x\left( \int_{]0,t]\times \{y\neq 0\}}  i \xi e^{i (X_{s-}-x)'\xi}y(1-\chi(y)) \ \mu^{X^\sigma}(\cdot;ds,dy)\right) \\
&+ \bbe^x \left(\int_{0+}^t i \xi e^{i(X^\sigma_{s-}-x)'\xi} \ dB^\sigma_s \right)\\
\end{align*}
First, we show that
\[
\bbe^x \left[\int_{0+}^t i \xi e^{i(X_{s-}-x)'\xi} \ d\left(X_s^c \right)^\sigma \right]= 0.
\]
The integral $e^{i(X_{t-}-x)'\xi}  \Si  X_t^c$ is a local martingale, since $X_t^c$ is a local martingale. To see that it is indeed a martingale, we calculate the following:
In the first two lines the integrand is now bounded, because $\ell$ and $Q$ are locally bounded and $\norm{X^{\sigma}_s(\omega)}<k$ on $\roi{0}{\sigma(\omega)}$ for every $\omega \in \Omega$.  
For the martingale preservation in the first term we obtain
\begin{align*}
\left[e^{i(X^\sigma-x)'\xi} \Si  X^{\sigma,c}, e^{i(X^\sigma-x)'\xi} \Si  X^{\sigma,c} \right]_t
&= \left[e^{i(X^\sigma-x)'\xi} \Si  X^{c,\sigma}, e^{i(X^\sigma-x)'\xi} \Si  X^{c,\sigma} \right]_t  \\
&= \left[(e^{i(X^\sigma-x)'\xi} \Si  X^{c})^\sigma, (e^{i(X^\sigma-x)'\xi} \Si  X^{c})^\sigma \right]_t \\
&= \left[e^{i(X^\sigma-x)'\xi} \Si  X^{c}, e^{i(X^\sigma-x)'\xi} \Si  X^{c} \right]^\sigma_t \\
&= \left(\int_0^t (e^{i(X^\sigma_s-x)'\xi})^2  \  d[X^{c},X^{c}]_s\right)^\sigma \\
&= \int_0^t \left(e^{i(X^\sigma_s-x)'\xi}\right)^2 1_{\ls 0 , \sigma \ls}(s) \ d[X^{c},X^{c}]_s \\
&=\int_0^t \left(e^{i(X^\sigma_s-x)'\xi}\right)^2 1_{\ls 0 , \sigma \ls}(s) \ d\left(\int_0^s Q(X_{r}) \ dF_r \right) \\
&= \int_0^t \left(e^{i(X^\sigma_s-x)'\xi}\right)^2 1_{\ls 0 , \sigma \ls}(s)Q(X_s)  \ dF_s
\end{align*}
where we used several well known facts about the square bracket.  Hence,
$$
\bbe^x \left(\left[e^{i(X^\sigma-x)'\xi} \Si  X^{\sigma,c}, e^{i(X^\sigma_t-x)'\xi} \Si  X_t^{\sigma,c} \right]_t\right) < \infty, \quad t\geq 0,
$$
{\cor which} provides that $e^{i(X^\sigma_t-x)'\xi} \Si  X_t^{\sigma,c}$ is an $L^2$-martingale which is zero at zero and therefore, its expected value is constantly zero.  \\
Moreover,  one easily sees that $e^{i(X_{s-}-x)'\xi}y\chi(y)$ is in the class $F_p^2$ of Ikeda and Watanabe {\cor (\cite{ikedawat}, Section 2.3)}, and we conclude that 
$$
\int_0^t\int_{y\neq 0} \Big(e^{i(X_{s-}-x)'\xi}  \chi(y)y \Big) \mu^{X^\sigma}(\cdot;ds,dy)-\nu^\sigma(\cdot;ds,dy))
$$
is also a martingale. Thus, in summation we obtain 
\begin{align*}
 \bbe^x \left( e^{i(X^\sigma_t-x)'\xi } -1 \right)
&=  \bbe^x \left(\int_{0+}^t i \xi e^{i(X^\sigma_{s-}-x)'\xi} \ dX^\sigma_s \right)+  \bbe^x\left(\frac{1}{2} \int_{0+}^t-\xi^2 e^{i(X^\sigma_{s-}-x)'\xi} \ dC^\sigma_s \right)\\
&+ \bbe^x\left( \int_{]0,t]\times \{y\neq 0\}}  \Big(e^{i (X_{s-}-x)'\xi}(e^{i\xi'y}-1 -i \xi'y\chi(y))\Big)\ \mu^{X^\sigma}(\cdot;ds,dy)\right) \\
&+ \bbe^x\left( \int_{]0,t]\times \{y\neq 0\}} \Big(  e^{i (X_{s-}-x)'\xi}( -i\xi'y \cdot(1-\chi(y)))\Big) \ \mu^{X^\sigma}(\cdot;ds,dy) \right)\\
&=  \bbe^x\left( \int_{]0,t]\times \{y\neq 0\}}  i \xi e^{i (X_{s-}-x)'\xi}y(1-\chi(y)) \ \mu^{X^\sigma}(\cdot;ds,dy)\right) \\
&+ \bbe^x \left(\int_{0+}^t i \xi e^{i(X^\sigma_{s-}-x)'\xi} \ dB^\sigma_s \right)+  \bbe^x\left(\frac{1}{2} \int_{0+}^t-\xi^2 e^{i(X^\sigma_{s-}-x)'\xi} \ dC^\sigma_s \right)\\
&+ \bbe^x\left( \int_{]0,t]\times \{y\neq 0\}}  \Big(e^{i (X_{s-}-x)'\xi}(e^{i\xi'y}-1 -i \xi'y\chi(y))\Big)\ \mu^{X^\sigma}(\cdot;ds,dy)\right) \\
&+ \bbe^x\left( \int_{]0,t]\times \{y\neq 0\}} \Big(  e^{i (X_{s-}-x)'\xi}( -i\xi'y \cdot(1-\chi(y)))\Big) \ \mu^{X^\sigma}(\cdot;ds,dy) \right)\\
&=\bbe^x \left(\int_{0+}^t i \xi e^{i(X^\sigma_{s-}-x)'\xi} \ dB^\sigma_s \right) +  \bbe^x\left(\frac{1}{2} \int_{0+}^t-\xi^2 e^{i(X^\sigma_{s-}-x)'\xi} \ dC^\sigma_s \right)\\
&+ \bbe^x\left( \int_{]0,t]\times \{y\neq 0\}}  \Big(e^{i (X_{s-}-x)'\xi}(e^{i\xi'y}-1 -i \xi'y\chi(y))\Big)\ \mu^{X^\sigma}(\cdot;ds,dy)\right)
\end{align*}
By the associativity of the stochastic integral and the properties of the compensator of a random measure the following equation holds true 
\begin{align*}
 \bbe^x \left( e^{i(X^\sigma_t-x)'\xi } -1 \right)
 &=\bbe^x \left[ \int_{0}^t \left(i \xi e^{i(X^\sigma_{s-}-x)'\xi}\ell(X_s)1_{\ls 0,\sigma \ls}-  \frac{1}{2}\xi^2 e^{i(X^\sigma_{s-}-x)'\xi} Q(X_s)1_{\ls 0,\sigma \ls} \right.\right.\\
&+ \left. \left.\int_{ \{y\neq 0\}}  \Big(e^{i (X_{s-}-x)'\xi}(e^{i\xi'y}-1 -i \xi'y\chi(y))\Big)1_{\ls 0,\sigma \ls} \ N(X_s(\omega),dy) \right)\ dF_s(\omega)\right].
\end{align*}
\end{proof}

{\cor
Before we state the next lemma we want to clear up one more notation: Let $(\mathcal{F}^0_t)_{t \geq 0}$ and $\mathcal{F}^0$ be defined as before. We define $(\mathcal{F}_t)_{t \geq 0}$ and $\mathcal{F}$ to be the completion of $(\mathcal{F}^0_t)_{t \geq 0}$ resp. $\mathcal{F}^0$ with respect to $\{ \bbp^\mu: \mu \text{ is a finite measure on }\bbr^d\}$ where $\bbp^\mu(\cdot) := \int \bbp^x(\cdot) \ \mu(dx)$. Note, that $(\mathcal{F}_t)_{t\geq 0}$ is right-continuous, and $\mathcal{F}\subset \mathcal{M}$, $\mathcal{F}_t \subset \mathcal{M}_t$.
\begin{lemma}\label{Lemma 2}
Let $(\Omega, \mathcal{M}, (\mathcal{M}_t)_{t \geq 0}, (X_t)_{t\geq 0}, (\theta_t)_{t \geq 0}, \bbp^x)_{x \in \bbr}$ be a Hunt process. Let $(F_t)_{t \geq 0}$ be an increasing, continuous and strongly AF adapted to $(\mathcal{F}_t)_{t\geq 0}$.\\
Then the processes $(F_t)_{t \geq 0}$ and $\left( \int_0^t g(X_s) ds + S_t \right)_{t \geq 0}$ are indistinguishable,
where $g$ is a positive $\mathcal{M}$-measurable function and $(S_t)_{t \geq 0}$ is an increasing, continuous, singular, perfect strongly AF adapted to $(\mathcal{F}_t)_{t \geq 0}$.
\end{lemma}
\begin{proof}
Let $B_t:=t$ for all $t \geq 0$.
Combination of Exercise 66.17 and Theorem 8.6 of \cite{sharpe} provides the existence of a set $G \in \mathcal{B}(\bbr^d)$ and a measurable function $g$ such that $1_G(X) \Si F = g(X) \Si B$. Moreover, the strongly additive functional $S:=1_{G^c}(X) \Si F$ is singular.  We obtain 
\begin{align*}
F = (1_G(X) + 1_{G^c}(X)) \Si F = g(X) \Si B + S.
\end{align*}
Since $F$ is increasing, continuous and adapted to $(\mathcal{F}_t)_{t\geq 0}$ so is $S$, and Proposition 3.21 of \cite{vierleute} provides that $S$ is indistinguishable from a perfect strongly additive functional.  
\end{proof}
}

\begin{theorem} \label{thm:main}
Let $(\Omega, \mathcal{M}, (\mathcal{M}_t)_{t \geq 0}, (X_t)_{t\geq 0}, (\theta_t)_{t \geq 0}, \bbp^x)_{x \in \bbr^d}$ be a Hunt semimartingale, i.e., a strongly Markovian, quasi-left-continuous $\bbp^x$-semimartingale.  Moreover, let the functions $\ell$ and $Q$ as mentioned in (\ref{b}) and (\ref{c})  be continuous, and let the function 
$$
x \mapsto \int_{\{y \neq 0\}} (1\wedge y^2) \ N(x,dy)
$$ be continuous.\\
If the symbol of $X$, namely 
\begin{equation}\label{1}
p(x,\xi) = - \lim_{t \downarrow 0} \bbe^x \left[ \frac{e^{i(X^\sigma_t-x)'\xi}-1}{t}\right]
\end{equation}
exists for all $x, \xi \in \bbr^d$, then $X$ is a Lévy-type process.
\end{theorem}

\begin{proof}
In this proof we only consider the one-dimensional case. The multivariate one works alike, but it is notationally more involved. \\
The stopping time $\sigma$ mentioned in (\ref{1}) is defined by 
$$
\sigma:= \inf \{ t \geq 0 : \| X_t-x \| > k \}
$$
for some starting point $x$ and $k \in \bbr_+$.  As we have already pointed out in the introduction, the limit does not depend on the particular choice of $k$. Hence, we will be able to choose a particular $k$ later in this proof. \\
At first, let us consider the expression 
$$
\bbe^x \left[ \frac{e^{i(X^\sigma_t-x)\xi}-1}{t}\right].
$$
Lemma $\ref{Lemma 1}$ provides the following equation:
\begin{align*}
\bbe^x \left[ e^{i(X^\sigma_t-x)\xi}-1\right]
=&\bbe^x \left[  \int_0^t i\xi e^{i(X_{s-}-x)\xi}\ell(X_s)1_{\ls 0, \sigma \ls}- \frac{1}{2}\xi^2 e^{i(X_{s-}-x)\xi} Q(X_s)1_{\ls 0, \sigma \ls}\right. \\
&\hspace{0.5cm} +\left. \int_{\{ y \neq 0 \}} e^{i(X_{s-}-x)\xi} 1_{\ls 0, \sigma \ls}\left(e^{i\xi y} -1 -i\xi y \chi(y)\right) \ N(X_s,dy)\ dF_s \right].
\end{align*}
So, when considering the symbol we obtain
\begin{align*}
p(x,\xi) &= - \lim_{t \downarrow 0} \frac{1}{t}  \bbe^x \left[ e^{i(X^\sigma_t-x)\xi}-1\right]\\
&= \lim_{t \downarrow 0} \bbe^x \left[  \int_0^t -i\xi e^{i(X_{s-}-x)\xi}\ell(X_s)1_{\ls 0, \sigma \ls}+ \frac{1}{2}\xi^2 e^{i(X_{s-}-x)\xi} Q(X_s)1_{\ls 0, \sigma \ls}\right. \\
&\hspace{0.5cm} -\left. \int_{\{ y \neq 0 \}} e^{i(X_{s-}-x)\xi} 1_{\ls 0, \sigma \ls}\left(e^{i\xi y} -1 -i\xi y \chi(y)\right) \ N(X_s,dy)\ dF_s \right].\\
\end{align*}
We define 
$$
Y_s^{\xi,x}:= -i\xi e^{i(X_{s}-x)\xi}\ell(X_{s})+ \frac{1}{2}\xi^2 e^{i(X_{s}-x)\xi} Q(X_{s})- \int_{\{ y \neq 0 \}} e^{i(X_{s}-x)\xi} \left(e^{i\xi y} -1 -i\xi y \chi(y)\right) \ N(X_s,dy). 
$$ 
Hence,  by the continuity of $F$ the existing limit in (\ref{1}) is given by
\begin{equation}\label{2}
p(x,\xi)= \lim_{t \downarrow 0} \bbe^x \left[ \int_0^t Y_s^{\xi,x} 1_{\ls 0, \sigma \ls} \ dF_s \right].
\end{equation}
Let us consider the process $(F_t)_{t \geq 0}$: by Lemma \ref{Lemma 2} {\cor this process is indistinguishable to 
$$
\left(\int_0^t g(X_u) du +S_t \right)_{t\geq 0}
$$
where $g$ is a positive, measurable function and $(S_t)_{t \geq 0}$ is an increasing, continuous,  singular, perfect strongly AF}. 
Thus, we decompose (\ref{2}) as follows: 
\begin{align}
p(x,\xi) = \lim_{ t\downarrow 0}\left( \frac{1}{t}\bbe^x \left[ \int_0^t Y_u^{\xi,x}g(X_u)1_{\ls 0, \sigma \ls} \ du \right]+ \frac{1}{t}\bbe^x \left[ \int_0^t Y_u^{\xi,x}1_{\ls 0, \sigma \ls} \ dS_u \right]\right) \label{3}
\end{align}
For the first summand of (\ref{3}), Theorem 4.4 of \cite{mydiss} provides the existence of the limit, and since (\ref{3}) exists, the limit 
\begin{equation*}
\lim_{ t\downarrow 0} \frac{1}{t}\bbe^x \left[ \int_0^t Y_u^{\xi,x}1_{\ls 0, \sigma \ls} \ dS_u \right] 
\end{equation*}
exists for all $\xi, x \in \bbr$, and, therefore, the same is true for
\begin{equation}\label{4}
\lim_{ t\downarrow 0} \frac{1}{t}\bbe^x \left[ \int_0^t \Re \left(Y_u^{\xi,x}\right)1_{\ls 0, \sigma \ls} \ dS_u \right],
\end{equation}
and 
\begin{equation}\label{4.2}
\lim_{ t\downarrow 0} \frac{1}{t}\bbe^x \left[ \int_0^t \Im \left(Y_u^{\xi,x}\right)1_{\ls 0, \sigma \ls} \ dS_u \right].
\end{equation}
First, we consider the case where there exists a $x_0 \in \bbr$ such that  $Q(x_0) \neq 0$ or $N(x_0, dy) \neq 0$:\\
the function $q^{\xi,x}$ defined by
\begin{align*}
z \mapsto  \Re \left( -i\xi e^{i(z-x)\xi}b(z)+\frac{1}{2}\xi^2 e^{i(z-x)\xi} c(z)-\int_\bbr e^{i(z-x)\xi} \left(e^{i\xi y} -1 -i\xi y \chi(y)\right) \ N(z,dy) \right)
\end{align*}
is continuous, and satisfies
$$
q^{\xi,x_0}(x_0)=\frac{1}{2}\xi^2 c(x_0)+\int_\bbr  1-\cos(y\xi)  \ N(x_0,dy) >0.
$$
The continuity provides the existence of  a $\delta>0$ such that $q^{\xi,x_0}(z) >\varepsilon$ for all $\| z- x_0 \| \leq \delta$ and an $\varepsilon >0$. As mentioned at the beginning of this proof we again consider the stopping time 
$$
\sigma:= \inf \{ t \geq 0 : \| X_t-x_0 \| > k \},
$$
and we set $k = \delta$. 
Thus, 
\begin{align*}
q^{\xi,x_0}(X_t)1_{\ls 0, \sigma \ls}> 0 \quad \text{ a.s.}
\end{align*}
Now, we again consider the limit in (\ref{4}):
$$
\lim_{ t\downarrow 0} \frac{1}{t}\bbe^{x_0} \left[ \int_0^t  \Re \left(Y_u^{\xi,x_0}\right)1_{\ls 0, \sigma \ls} \ dS_u \right]=\lim_{ t\downarrow 0} \frac{1}{t}\bbe^{x_0} \left[ \int_0^t  \Re \left(q^{\xi,x_0}(X_u) \right) 1_{\ls 0, \sigma \ls}\ dS_u \right]
$$
which exists for all $\xi,x\in \bbr$. Moreover $(S_t)_{t \geq 0}$ is a path-wise singular function.  
In order to establish the statement, we want to show that $(S_t)_{t \geq 0} \equiv 0$ (in this case $X$ is a Lévy-type process).  Thus, we assume that there exists a set $M \in \mathcal{M}$ with $\bbp(M) > 0$ such that $S_\cdot (\omega) \neq 0$ for all $\omega \in M$. \\
By Proposition \ref{prop:derived} for all $\omega \in M$ there exists a $t(\omega) \geq 0$ with 
$$
D_+ S_{t(\omega)} = \infty.
$$
We define 
$$
\tau := \inf\{ u \geq 0: D_+ S_u = \infty \} \quad \forall\omega \in \Omega
$$
which is an $( \mathcal{M}_t)_{t \geq 0}$-stopping time since $D_+S_u$ is $ \mathcal{M}_{u+}$-measurable and the filtration is right-continuous.  
Thus, 
\begin{equation}\label{5}
\lim_{t \downarrow 0} \frac{1}{t} \bbe^{X_{\tau(\omega)}(\omega)} \left[ \int_0^t \Re \left(Y_u^{\xi,x_0}\right)1_{\ls 0, \sigma \ls}  \ dS_u\right]
\end{equation}
exists for all $x,\xi \in \bbr$ and $\omega \in M$. \\ 
By applying Fatou's lemma for the conditional expectation the following inequality holds true 
\begin{align*}
&\lim_{t \downarrow 0} \frac{1}{t} \bbe^{X_{\tau(\omega)}(\omega)} \left[ \int_0^t \Re \left(Y_u^{\xi,x_0}\right)1_{\ls 0, \sigma \ls}  \ dS_u\right] \\
=&\lim_{t \downarrow 0} \frac{1}{t} \bbe^{x_0} \left[ \int_0^t \Re \left(Y_u^{\xi,x_0}\circ \theta_\tau \right)1_{\ls 0, \sigma \ls}\circ \theta_\tau  \ dS_u\circ \theta_\tau \; \bigg| \; \mathcal{M}_\tau\right] \\
\geq &\bbe^{ x_0} \left[ \liminf_{t \downarrow 0} \frac{1}{t}\int_0^t \Re \left(Y_u^{\xi,x_0}\circ \theta_\tau \right)1_{\ls 0, \sigma \ls}\circ \theta_\tau  \ dS_u\circ \theta_\tau \; \bigg| \; \mathcal{M}_\tau\right] \\
\geq &\bbe^{ x_0} \left[ \liminf_{t \downarrow 0} \left( \inf_{u \in [0,t]} \Re \left(Y_u^{\xi,x_0}\circ \theta_\tau \right)1_{\ls 0, \sigma \ls}\circ \theta_\tau  \right) \ \frac{S_t\circ \theta_\tau}{t} \; \bigg| \; \mathcal{M}_\tau \right] \\
= &\bbe^{ x_0} \left[ \liminf_{t \downarrow 0} \left( \inf_{u \in [0,t]} \Re \left(Y_u^{\xi,x_0}\circ \theta_\tau \right)1_{\ls 0, \sigma \ls}\circ \theta_\tau  \right) \ \frac{S_{\tau + t} - S_\tau }{t} \; \bigg| \; \mathcal{M}_\tau \right] 
\end{align*}
{\cor where we have used the strong additivity and the perfectness of $S$ in the last equation}.\\
Since $\Re \left(Y_u^{\xi,x_0}\circ \theta_\tau \right)1_{\ls 0, \sigma \ls}\circ \theta_\tau  >\varepsilon$ and $(S_{\tau + t} - S_\tau)/t \to \infty $ for $t\downarrow 0$ by the definition of $\tau$ this is a contradiction.\\
In case that $c, N \equiv 0$ there exists  $x_0 \in \bbr$ with $b(x_0) \neq 0$. Otherwise the processes would be constantly zero. We set
\begin{align*}
q^{\xi,x}(z)&:= \Im \left( i\xi e^{i(z-x)\xi}b(z)- \frac{1}{2}\xi^2 e^{i(z-x)\xi} c(z)-\int_\bbr e^{i(z-x)\xi} \left(e^{i\xi y} -1 -i\xi y \chi(y)\right) \ N(z,dy) \right)\\
&= \xi \cos \left((z-x)\xi\right)b(z)
\end{align*}
which is continuous, and satisfies
$$
q^{\xi_0,x_0}(x_0)=\xi_0 b(x_0) >0 \quad a.s.
$$
for $\xi_0 = \sign \left(  b(x_0) \right)$.  
The rest of the proof is analogous to the case above but using of (\ref{4.2}) instead of (\ref{4}). 
\end{proof}

\section{Generalizations Beyond Hunt Processes} \label{sec:final}

As we have pointed out above it was believed that for the symbol to exist in the Markovian framework, it was necessary to consider a L\'evy-type process. Now we construct a counter-example, that is, a Markov semimartingale, which is not L\'evy-type but having a symbol.
Most people would probably start with Brownian Motion, transforming and/or time-shifting it. Our process is different, since it does not even depend on random, i.e.,  it is a deterministic function. Rather one should say `a family of functions' satisfying the deterministic version of \eqref{universalmp}. We have dealt with processes of this kind in \cite{detmp2} and \cite{detmp1}. In the following we will work one-dimensional.

\begin{lemma}
Given a family of functions $f_x:[0,\infty[\to\bbr$, $x\in\bbr$, the corresponding stochastic process $X_t^x:=f_x(t)$ is a time-homogeneous Markov process, iff for every $s,t\geq 0$ and $x,y\in\bbr$ such that $f_x(s)=f_y(t)$ we obtain
\begin{align} \label{timehom}
  f_x(s+h)=f_y(t+h)
\end{align}
for $h\geq 0$.
\end{lemma}

For a deterministic process - which is allowed to start in every point $x\in\bbr$ - the symbol boils down to 
\[
-\lim_{t\downarrow 0} \left(\frac{\cos((X^x_t-x)\xi)-1}{t}+ i\frac{\sin((X^x_t-x)\xi)}{t}\right).
\]
Whether we use the version \eqref{probsymbol1} or \eqref{probsymbol2} does not make any difference here. Furthermore, since (for $\xi\neq 0$) the function $y\mapsto \sin((y-x)\xi)$ is locally {\cor diffeomorphic} around zero, the symbol exists if and only if the process is right-differentiable at zero in every starting point $x$. In this case the cosine-part becomes zero. 

\begin{example} \label{example:det1}
Writing as usual every $y\in\bbr$ as $y=[y]+\{y\}$, where $[\cdot]$ denotes the floor function and $\{y\}\in[0,1[$ we define
\begin{align}\label{loop}
X_t^x:= [x]+\{x+t\}.
\end{align}
Starting from $x$ the process drifts with constant speed upwards, but jumps to $[x]$ directly before reaching $[x]+1$. The symbol of this process exists and is equal to $-i\xi$. This means that {\cor the symbol of this process is exactly the one of a deterministic drift, that is, the most simple example of a L\'evy process.}

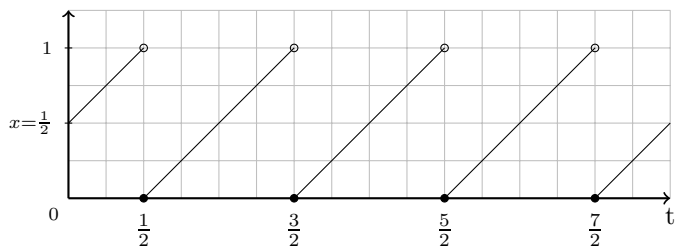
\begin{center}
\begin{tikzpicture}[x=.5cm, y=.5cm,domain=-9:9,smooth]
   \draw [color=gray!50]  [step=5mm] (0,0) grid (16,5);
   \draw[->,thick] (0,0) -- (16,0) node[below]{t};
   \draw[->,thick] (0,0) -- (0,5) node[above]{} ;
  
     \draw (2,-.1) -- (2,.1) node[below=4pt] {$\frac{1}{2}$};
     \draw (6,-.1) -- (6,.1) node[below=4pt] {$\frac{3}{2}$};
      \draw (10,-.1) -- (10,.1) node[below=4pt] {$\frac{5}{2}$};
      \draw (14,-.1) -- (14,.1) node[below=4pt] {$\frac{7}{2}$};

     \draw (-.1,2) -- (.1,2) node[left=4pt] {$\scriptstyle x=\frac{1}{2}$};
      \draw (-.1,4) -- (.1,4) node[left=4pt] {$\scriptstyle 1$};
   
   \node[below left]{$\scriptstyle0$};
\draw (0,2) -- (2,4);
\draw (2,0) -- (6,4);
\draw (6,0) -- (10,4);
\draw (10,0) -- (14,4);
\draw (14,0) -- (16,2);
\draw[fill=black] (2,0) circle(0.5mm) ;
\draw[fill=black] (6,0) circle(0.5mm) ;
\draw[fill=black] (10,0) circle(0.5mm) ;
\draw[fill=black] (14,0) circle(0.5mm) ;
\draw (2,4) circle(0.5mm) ;
\draw(6,4) circle(0.5mm) ;
\draw (10,4) circle(0.5mm) ;
\draw (14,4) circle(0.5mm) ;
    
\end{tikzpicture} 
\end{center}
Admittedly, this process does (at first) not look Markovian, since in the definition \eqref{loop} the starting point $x$ appears. In fact this is just to simplify notation. Being at time $s$ at the point $y$, we know that at time $t>s$ we are at $[y]+\{y+(t-s)\}$, without any knowledge on the fact where the process has started at time zero.
\end{example}

We have seen above that under some mild technical assumptions within the Hunt framework the following holds: the process admits a probabilistic symbol iff it is a L\'evy-type process.  However, in {\cor Example \ref{example:det1}} we have constructed a Markov semimartingale, which is not L\'evy-type but which admits a symbol.  {\cor We will now see}, that the symbol does not contain the same information on the process as before. We focus on one particular property of the paths, but the same holds true for various other properties of the process. 

\begin{example} \label{example:det2}

Let us for simplicity consider again a deterministic Markov process as in Example \ref{example:det1}.  Starting at time zero in zero, the paths looks at follows:
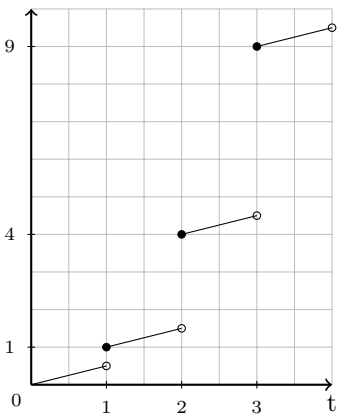
\begin{center}
\begin{tikzpicture}[x=.5cm, y=.5cm,domain=-9:9,smooth]
   \draw [color=gray!50]  [step=5mm] (0,0) grid (8,10);
   \draw[->,thick] (0,0) -- (8,0) node[below]{t};
   \draw[->,thick] (0,0) -- (0,10) node[above]{} ;

     \draw (2,-.1) -- (2,.1) node[below=4pt] {$\scriptstyle 1$};
      \draw (4,-.1) -- (4,.1) node[below=4pt] {$\scriptstyle 2$};
     \draw (6,-.1) -- (6,.1) node[below=4pt] {$\scriptstyle 3$};
      
   \foreach \c in {1,4,9}{
     \draw (-.1,\c) -- (.1,\c) node[left=4pt] {$\scriptstyle\c$};
     }
     
   \node[below left]{$\scriptstyle0$};
\draw (0,0) -- (2,0.5);
\draw (2,1) -- (4,1.5);
\draw (4,4) -- (6,4.5);
\draw (6,9) -- (8,9.5);
\draw[fill=black] (2,1) circle(0.5mm) ;
\draw[fill=black] (4,4) circle(0.5mm) ;
\draw[fill=black] (6,9) circle(0.5mm) ;
\draw (2,0.5) circle(0.5mm) ;
\draw(4,1.5) circle(0.5mm) ;
\draw (6,4.5) circle(0.5mm) ;
\draw (8,9.5) circle(0.5mm) ;
    
\end{tikzpicture} 
\end{center}
A drift with incline 1/2 is disrupted by deterministic jumps. The paths from other starting points are added in a time-homogeneous manner. 
Starting in zero, for $t\to \infty$ our process behaves like $t^2$. 
The symbol of the process is: $p(x,\xi)=\psi(\xi)=-\ell\xi$, where in this special case $\ell=(1/2)$. 
In the class of L\'evy-type processes one considers 
\[
H(R):=\sup_{y\in\bbr} \sup_{\norm{\varepsilon}\leq 1} \abs{p\left( y, \frac{\varepsilon}{R} \right)} = \frac{\ell}{R},
\]
and derives the so called generalized Blumenthal-Getoor index as follows: 
\[
\beta_0:=\sup\left\{\lambda \geq 0: \limsup_{R\to\infty} R^\lambda H(R)=0\right\} = 1.
\]
A typical application for indices of this kind deals with the behavior of the paths when $t$ tends to infinity.  For L\'evy-type processes it holds that (cf. \cite{generalizedindizes} Theorem 3.11)
\[
\lim_{t\to\infty} t^{-1/\lambda} (X-x)_t^*=0 \text{ for all } \lambda<\beta_0.
\]
This is obviously wrong for the process considered here, although the symbol and the corresponding index do exist.  In fact, one can mimic totally different kinds of behavior for $t\to\infty$, by using different jump structures.  Various other results, like the maximal inequality \cite{generalizedindizes} Proposition 3.10 fail to hold also. 
\end{example}

If one needs a probabilistic example (for some reason), one can make e.g. the jump sizes random.  In any case, the example shows that if we leave {\cor the Hunt framework} behind, the symbol might exist, but is of no use anymore.  Hence, within the framework where the symbol is useful, we have found a characterization for its existence. \vspace{5mm}

{\cor \textbf{Acknowledgement}: We would like to thank an anonymous referee whose comments significantly helped us to improve the paper. }

\appendix

\section{Appendix}

Let us first recall the following definition:

\begin{definition}\label{def:derivednr}
The number $c$ (finite or infinite) is called a \emph{derived number} of the function $f$ at the point $x_0$ if there exists a sequence $h_n\to 0$ and 
\[
\lim_{n\to \infty} \frac{f(x_0+h_n)-f(x_0)}{h_n}=c.
\]
\end{definition}

The {\cor proof of Proposition \ref{prop:derived}} uses two lemmas which we recall subsequently. We write $\lambda$ for Lebesgue measure {\cor defined on the Lebesgue sets, that is, we are working on the completion of the Borel sets w.r.t. Lebesgue measure}. $\lambda^*$ denotes the corresponding outer measure. 
In fact, all sets on which we will use the following results will be measurable.   The first result is  \cite{kankrue}  Lemma 1.2.3 while the second one is a combination of Lemma 1.2.3 and Lemma 1.2.5 of \cite{kankrue}.

\begin{lemma} \label{lem:lemma2}
Let $f$ be a strictly increasing function on $[a,b]$ and let $p\geq 0$. If at every point $x$ of a set $E\subseteq [a,b]$, there exists at least one derived number $Df(x)$ such that $Df(x)\leq p$ then
\[
\lambda^*(f(E))\leq p \cdot \lambda^*(E).
\]
\end{lemma}

\begin{lemma} \label{lem:lemma3}
Let $f$ be a strictly increasing function. If at every point $x$ of the set $E\subseteq [a,b]$, there exists $f'(x)=p$ then
\[
\lambda^*(f(E))=p\cdot \lambda^*(E).
\]
\end{lemma}

\begin{proof}[Proof of Proposition \ref{prop:derived}]

We prove the result only for the upper right Dini derivative $D^+$ since all other cases, i.e. , $D_+,D^-,D_-$, work analogously.\\
Let $t \in \bbr_+$ with $g(t) \neq 0$. 
We consider $g$ on the interval $I:=[t-0.5, t+0.5 ]$. Let us assume that $g$ possesses a finite upper right Dini derivative in every point of $I$. For $x \in I$, we consider $f(x):=g(x)+x$ which is a strictly increasing, continuous function with $f'(x)=1$ almost everywhere, i.e., in a set of measure $1$. We define the disjoint sets
\begin{align*}
B&:= \{ x \in I : f'(x)=1\},\\
E_j &:= \{ x \in I : j \leq D^+f(x) < j+1\} \text{ for } j \in \bbn,
\end{align*}
and obtain that $I = B \cup \left( \bigcup_{j=1}^\infty E_j \right)$.  Moreover, the sets $B$ and $(E_j)_{j \in \bbn}$ are measurable by Theorem 3.6.5 of \cite{kankrue}, and since $f$ is continuous and strictly increasing the sets $f(B),f(E_1),...$  form a disjoint decomposition of $f(I)$ into measurable sets. The sets $(E_j)_{j \in \bbn}$ do all have Lebesgue measure zero.  We obtain 
\begin{align*}
\lambda(f(I)) &= \lambda(f(B  \cup E_1 \cup...)) \\
&=\lambda(f(B)  \cup f(E_1) \cup...) \\
&= \lambda(f(B))\\
&= \lambda (B)\\
&=1 
\end{align*}
where we used Lemma \ref{lem:lemma2} in the third and Lemma \ref{lem:lemma3} in the fourth equation.  Thus, since $f$ maps intervals to intervals, we obtain that  $f(t-0.5)=c$ and $f(t+0.5)=c+1$ for a $c \in \bbr$,  and conclude that $g \equiv 0$ on $I$ by monotonicity. This is a contradiction.
\end{proof}

\end{document}